\documentclass[a4paper,12pt]{amsart}
\usepackage{ifthen}
\usepackage[dvipdfmx]{graphicx}
\usepackage{mathrsfs}
\nonstopmode \numberwithin{equation}{section}
\newtheorem{thm}{Theorem}[section]
\newtheorem{cor}[thm]{Corollary}
\newtheorem{lem}[thm]{Lemma}

\newtheorem{Cor}{Corollary}

\theoremstyle{definition}

\newtheorem{example}[thm]{Example}
\newtheorem{rem}[thm]{Remark}

\newenvironment{pf}[1][]{%
 \vskip 3mm
 \noindent
 \ifthenelse{\equal{#1}{}}%
  {{\slshape Proof. }}%
  {{\slshape #1.} }%
 }%
{\qed\bigskip}

\newcounter{alphabet}
\newcounter{tmp}
\newenvironment{Thm}[1][]{\refstepcounter{alphabet}%
\bigskip%
\noindent%
{\bf Theorem \Alph{alphabet}}%
\ifthenelse{\equal{#1}{}}{}{ (#1)}%
{\bf .} \itshape}{\vskip 8pt}

\newcommand{\A}{{\mathcal A}}

\newcommand{\C}{{\mathbb C}}
\newcommand{\D}{{\mathbb D}}

\newcommand{\R}{{\mathbb R}}

\newcommand{\es}{{\mathcal S}}

\newcommand{\spl}{{\mathcal{SP}}}

\newcommand{\st}{{\mathcal{SS}}}

\renewcommand{\Im}{{\,\operatorname{Im}\,}}

\renewcommand{\Re}{{\,\operatorname{Re}\,}}

\newcommand{\inv}{^{-1}}

\newcommand{\Gauss}{{\null_2F_1}}

\renewcommand{\arg}{\,{\operatorname{arg}\,}}

\newcommand{\aand}{{\quad\text{and}\quad}}
\newcommand{\sgn}{{\operatorname{sgn}}}

\newcounter{minutes}
\divide\time by 60
\newcounter{hours}
\multiply\time by 60 \addtocounter{minutes}{-\time}

\begin{document}
\bibliographystyle{amsplain}
\title{Spirallikeness of shifted hypergeometric functions}

\author[T. Sugawa]{Toshiyuki Sugawa}
\address{Graduate School of Information Sciences, \\
Tohoku University, \\
Aoba-ku, Sendai 980-8579, Japan} \email{sugawa@math.is.tohoku.ac.jp}
\author[L.-M.~Wang]{Li-Mei Wang}
\address{School of Statistics,
University of International Business and Economics, No.~10, Huixin
Dongjie, Chaoyang District, Beijing 100029, China}
\email{wangmabel@163.com}

\begin{abstract}
In the present paper, we study spirallikenss (including starlikeness)
of the shifted hypergeometric function $f(z)=z\Gauss(a,b;c;z)$ with
complex parameters $a,b,c,$ where
$\Gauss(a,b;c;z)$ stands for the Gaussian hypergeometric function.
First, we observe the asymptotic behaviour of $\Gauss(a,b;c;z)$ around the
point $z=1$ to obtain necessary conditions for $f$ to
be $\lambda$-spirallike for a given $\lambda$ with $- \pi/2<
\lambda<\pi/2.$
We next give sufficient conditions for $f$ to be $\lambda$-spirallike.
As special cases, we obtain sufficient conditions of strong starlikeness
and examples of spirallike, but not starlike, shifted hypergeometric functions.
\end{abstract}
\thanks{
} \maketitle

\section{Introduction and main results}
The {\it Gaussian hypergeometric function} $\Gauss(a,b;c;z)$
with complex parameters $a,b,c$~ $(c\not=0,-1,-2,\dots)$
is defined by the power series
$$
\Gauss(a,b;c;z)=\sum_{n=0}^{\infty}\frac{(a)_n(b)_n}{(c)_nn!}z^n
$$
for $z\in\D=\{z\in\C: |z|<1\},$
where $(a)_n$ is the {\it Pochhammer symbol}; namely, $(a)_0=1$ and
$(a)_n=a(a+1)\cdots(a+n-1)=\Gamma(a+n)/\Gamma(a)$ for $n=1,2,\dots.$
It is well known that $\Gauss(a,b;c;z)$ analytically extends to
the slit plane $\C\setminus[1,+\infty).$
For basic properties of hypergeometric functions, one can consult
\cite{AS:hand}, \cite{Temme:sp} or \cite{WW:anal}.

Let $\A$ denote the set of analytic functions $f$ on the
open unit disk $\D$
and consider the subclass $\A_1=\{f\in\A: f(0)=f'(0)-1=0\}.$
We denote by $\es$ the subset of $\A_1$ consisting of
univalent functions on $\D.$
For a constant $\lambda\in(-\pi/2,\pi/2)$, a function $f\in\A_1$ is
called {\it $\lambda$-spirallike} if
$$
\Re\left(e^{-i\lambda}\frac{zf'(z)}{f(z)}\right)>0,\quad z\in\D.
$$
(Note that in the literature a $\lambda$-spirallike function may refer
to $(-\lambda)$-spirallike one in our definition.)
Let $\spl(\lambda)$ denote the class of
$\lambda$-spirallike functions.
It is known that $\spl(\lambda)\subset\es.$
For a geometric characterization and other properties of
$\lambda$-spirallike functions, the reader may refer to \cite{Duren:univ}
(and also \cite{KS12spl}).
In particular, a function in
$\spl(0)$ is called {\it starlike} and we sometimes write
$\es^{*}=\spl(0).$
Furthermore, let
$$
\sigma(f)=\inf_{z\in\D}\Re\left(\frac{zf'(z)}{f(z)}\right).
$$
Here, following the convention adopted by K\"ustner \cite{Kus07},
we will leave $\sigma(f)$ undefined if $zf'(z)/f(z)$ has a pole in $\D$
so that the assertion $\sigma(f)=-\infty$ means that $zf'(z)/f(z)$
is pole-free but its real part has no lower bound on $\D.$
A function $f\in\A_1$ is called {\it starlike of order $\alpha$}
if $\sigma(f)\ge\alpha.$
Note here that $f$ is starlike precisely if $\sigma(f)\ge0.$
A (not necessarily normalized) function $f\in\A$ is called {\it
convex} if $f$ maps $\D$ univalently onto a convex domain.
It is well known that $f$ is convex if and only if $\Re[1+zf''(z)/f'(z)]>0$
on $|z|<1.$
For a real constant $\alpha\in(0,1)$, a function
$f\in\A_1$ is called {\it strongly starlike of order $\alpha$} if
$$
\left|\arg\frac{zf'(z)}{f(z)}\right|<\frac{\pi}{2}\alpha,\quad
z\in\D.
$$
Note that a strongly starlike function is starlike
and known to have a quasiconformal extension to the whole plane.
We denote by $\st(\alpha)$ the set of strongly starlike functions
of order $\alpha.$
For geometric properties of strongly starlike functions, the reader
may refer to \cite{SugawaDual} and cited papers there.
We note that for $\lambda\in(-\pi/2,\pi/2),$
a function $f\in\A_1$ is $\lambda$-spirallike if and only if
$$
\lambda-\frac\pi2<\arg\frac{zf'(z)}{f(z)}<\lambda+\frac\pi2.
$$
In particular, we observe that
$$
\spl(\lambda)\cap\spl(-\lambda)=\st(\alpha),\quad
\alpha=1-\frac2\pi|\lambda|.
$$

Note that the function $z\Gauss(a,b;c;z),$ called the {\it shifted
hypergeometric function}, belongs to the class $\A_1.$
A number of authors have investigated geometric properties of the shifted
hypergeometric functions.
For instance, sufficient conditions for those functions to be starlike or convex
were found by
Merkes and Scott \cite{MS61},
Lewis \cite{Lewis79}, Ruscheweyh and Singh \cite{RS86},
Miller and Mocanu \cite{MM90U},
Silverman \cite{Sil93},
Ponnusamy and Vuorinen \cite{PV01},
K\"ustner \cite{Kus02}, \cite{Kus07},
H\"asto, Ponnusamy and Vuorinen \cite{HPV10}.
Most of known results in this line, however, deal with
$z\Gauss(a,b;c;z)$ for real parameters $a,b,c$ only.
A few exceptions are \cite[Theorem 2.12]{Rus:conv}
(see also \cite[Theorem 4]{Kus07}), \cite[Theorem 14, Corollary 17]{Kus07}
(and its convex counterparts), and \cite[Remark 1.5]{HPV10}.
Moreover, to the best knowledge of the authors,
no results are found on spirallikeness of hypergeometric functions.
Some of known starlikeness results are summarized in the following form.
Note here that the hypergeometric functions are symmetric
in regard of the parameters $a$ and $b;$ namely,
$\Gauss(a,b;c;z)=\Gauss(b,a;c;z).$

\begin{Thm}
Let $f(z)=z\Gauss(a,b;c;z).$
Then the following hold:
\renewcommand{\labelenumi}{({\roman{enumi}})}
\begin{enumerate}
\item $($K\"ustner \cite[Remark 2.3]{Kus02}$)$
For $a,b,c\in\R$ with $0<a\le b\le c,$
$$
\sigma(f)=1-\frac{\Gauss'(a,b;c;-1)}{\Gauss(a,b;c;-1)}\ge 1-\frac{ab}{b+c}.
$$
\item $($K\"ustner \cite[Remark 1.2]{Kus02}$)$
For $a,b,c\in\R$ with $-1\le a<0<b$ and $c-a-b>1,$
$$
\sigma(f)=1-\frac{ab}{c-a-b-1}.
$$
\item $($Ruscheweyh, cf.~\cite[Theorem 4]{Kus07}$)$
For $a\in\R, b,c\in\C$ with $2\Re b\le a+1,~0\le a$ and $c=a-\bar b+1,$
$$
\sigma(f)\ge 1-\frac a2.
$$
\end{enumerate}
\end{Thm}

In particular, we obtain the following.

\begin{Cor}
The shifted hypergeometric function $z\Gauss(a,b;c;z)$ is starlike
if one of the following conditions is satisfied:
\renewcommand{\labelenumi}{({\roman{enumi}})}
\begin{enumerate}
\item 
$a,b,c\in\R$ with $0<a\le b\le c$ and $ab\le b+c.$
\item 
$a,b,c\in\R$ with $-1\le a<0<b$ and $c-a-b\ge 1-ab~ (>1).$
\item 
$a\in\R$ with $2\Re b\le a+1,~0\le a\le 2$ and $c=a-\bar b+1.$
\end{enumerate}
\end{Cor}

We remark that part (ii) in the corollary was first proved by
Silverman \cite{Sil93}.
In the present note, we study spirallikenss, including starlikeness,
of shifted hypergeometric functions with complex parameters.
First, we collect necessary conditions for spirallikeness
by looking at the behaviour as $z\to1$ in $\D.$
Since $f(z)\equiv z$ when $ab=0,$ it is reasonable to assume $ab\ne0$
from the beginning.

\begin{thm}\label{thm:nec}
We set $f(z)=z\Gauss(a,b;c;z)$ for complex numbers $a,b,c$ with
$ab\ne0, c\ne0,~ -1,$ $-2, \dots$ and let $-\pi/2<\lambda<\pi/2.$
Suppose that $f$ is $\lambda$-spirallike.
Then the following hold according to the value of $c-a-b:$
\renewcommand{\labelenumi}{({\roman{enumi}})}
\begin{enumerate}
\item
If $\Re(c-a-b)>1,$ then
$$
\left|\lambda-\arg\left(1+\frac{ab}{c-a-b-1}\right)\right|\le\frac\pi2.
$$
\item
If $c-a-b=1+si$ with $s\in\R\setminus\{0\},$ then $R_1\le |w_1|$ and
$$
\left|\lambda-\arg w_1\right|\le \arccos\frac{R_1}{|w_1|},
$$
where $w_1=1-iab/s$ and
$$
R_1
=\left|\frac{\Gamma(c-a)\Gamma(c-b)}{s\Gamma(a)\Gamma(b)}\right|e^{\pi|s|/2}
=\left|\frac{(a+is)(b+is)}s\cdot
\frac{\Gamma(a+is)\Gamma(b+is)}{\Gamma(a)\Gamma(b)}\right|e^{\pi|s|/2}.
$$
\item
If $c-a-b=1,$ then
$$
|\lambda-\arg(ab)|\le\frac\pi2.
$$
\item
If $0\le\Re(c-a-b)<1,$ then $c-a-b\in\R$ and
$$
\left|\lambda-\arg
\frac{\Gamma(c-a)\Gamma(c-b)}{\Gamma(a)\Gamma(b)}
\right|\le (c-a-b)\frac\pi2.
$$
\item
If $\Re(c-a-b)<0,$
$$
\lambda=\arg(a+b-c).
$$
\end{enumerate}
\end{thm}

We remark that the condition $R_1\le |w_1|$ in case (ii) is indeed
necessary for local univalence of the function $f(z).$
We now give several sufficient conditions for $f$ to be univalent.
The first result compliments Theorem A by adding a case of complex parameters.

\begin{thm}\label{thm:st}
Assume that complex numbers $a,b,c$ satisfy
$ab\ne0$ and $c\ne0,-1,-2,\dots.$
The shifted hypergeometric function $z\Gauss(a,b;c;z)$ is starlike
if the following conditions are satisfied:
\renewcommand{\labelenumi}{({\roman{enumi}})}
\begin{enumerate}
\item
$p=a+b+1-c$ is a real number,
\item
$\Re[ab]> p,$
\item
$L\ge0, N\ge0$ and $LN-M^2\ge0,$ where
\begin{align*}
L&=|c-1|^2-|a+b|^2+p+3\Re[ab], \\
M&=\Im[ab(\bar a+\bar b-2)], \aand \\
N&=|c-2|^2-|a-1|^2|b-1|^2-p+\Re[ab].
\end{align*}
\end{enumerate}
\end{thm}

The following simple fact might be helpful to check condition (iii).
The condition $N\ge0$ follows from the two inequalities
$L>0$ and $LN-M^2\ge0$ because $N\ge M^2/L\ge0.$
The roles of $L$ and $N$ are interchangeable.

By using Alexander's correspondence (see Lemma \ref{lem:Alex} given below),
the starlikeness criterion can readily be translated into a convexity one.


\begin{thm}\label{thm:convex}
Assume that complex numbers $a,b,c$ satisfy
$ab\ne0$ and $c\ne0,-1,-2,\dots.$
The hypergeometric function $\Gauss(a,b;c;z)$ is convex
if the following conditions are satisfied:
\renewcommand{\labelenumi}{({\roman{enumi}})}
\begin{enumerate}
\item
$p=a+b-c$ is a real number,
\item
$\Re[(a-1)(b-1)]> p,$
\item
$L\ge0, N\ge0$ and $LN-M^2\ge0,$ where
\begin{align*}
L&=|c-2|^2-|a+b-2|^2+p+3\Re[(a-1)(b-1)], \\
M&=\Im[(a-1)(b-1)(\bar a+\bar b-4)], \aand \\
N&=|c-3|^2-|a-2|^2|b-2|^2-p+\Re[(a-1)(b-1)].
\end{align*}
\end{enumerate}
\end{thm}

One might expect that the condition (ii) in Theorems \ref{thm:st} 
and \ref{thm:convex} could be weakened to allow equality.
This is indeed possible to some extent but not in full generality.
See Remarks \ref{rem:eq} and \ref{rem:ex} below.

A sufficient condition for spirallikeness can also be given as follows.
As long as we apply Jack's lemma in the present setting, it seems inevitable
to assume the additional condition $c=a+b+1$ (see the proof given in Section 3).

\begin{thm}\label{thm:spl}
Let $\lambda$ be a real number with $0<|\lambda|<\pi/2$
and $a,b$ be complex numbers.
Then the shifted hypergeometric function $z\Gauss(a,b;a+b+1;z)$ is
$\lambda$-spirallike if the following conditions are satisfied:
\renewcommand{\labelenumi}{({\roman{enumi}})}
\begin{enumerate}
\item
$\Re\big[e^{-i\lambda}ab\big]\ge0,$
\item
$L\ge0, N\ge0$ and $LN-M^2\ge0,$ where
\begin{align*}
L&=\Re\big[e^{-i\lambda}ab(2+e^{-2i\lambda})\big], \\
M&=\Im\big[e^{-i\lambda}ab(\bar a+\bar b-2e^{-i\lambda}\cos\lambda)\big],
\aand \\
N&=\Re\big[e^{-i\lambda}ab(2\bar a+2\bar b-e^{-2i\lambda}
-e^{i\lambda}\bar a\bar b/\cos\lambda)\big].
\end{align*}
\end{enumerate}
\end{thm}

We note that the function $f$ under the assumptions in Theorem \ref{thm:spl}
is always bounded (see Lemma \ref{lem:bdd}).
When $e^{-i\lambda}ab$ or $ab$ is real, the conditions in the theorem
may be simplified as follows.

\begin{cor}\label{cor:1}
Let $\lambda$ be a real number with $0<|\lambda|<\pi/2.$
Suppose that $q=e^{-i\lambda}ab$ is a positive real number.
Then the shifted hypergeometric function $z\Gauss(a,b;a+b+1;z)$
is $\lambda$-spirallike if
\begin{equation*}
(2+\cos2\lambda)(2\Re[a+b]-\cos2\lambda-q/\cos\lambda)
-(\Im[a+b]-\sin2\lambda)^2\ge0.
\end{equation*}
\end{cor}

\begin{cor}\label{cor:2}
Let $\lambda$ be a real number with $0<|\lambda|<\pi/3.$
Suppose that $q=ab$ is a positive real number.
Then the shifted hypergeometric function $z\Gauss(a,b;a+b+1;z)$
is $\lambda$-spirallike if
\begin{align}\label{eq:spl}
&(\Im[e^{i\lambda}(a+b)]-2\sin2\lambda\cos\lambda)^2 \\
\le~
&(4\cos^2\lambda-1)(2\Re[e^{i\lambda}(a+b)]\cos\lambda-4\cos^4\lambda
+3\cos^2\lambda-q)
\notag
\end{align}
\end{cor}

We finally obtain a sufficient condition for strong starlikeness.

\begin{thm}\label{thm:ss}
Let $1/3<\alpha<1$ and $a,b$ be complex numbers with
$a+b\in\R$ and $ab>0.$
Then the shifted
hypergeometric function $z\Gauss(a,b;a+b+1;z)$ is strongly starlike
of order $\alpha$ if
\begin{align}\label{eq:ss}
&~
[(a-b)^2+6(a+b)-3]\sin^2\frac{\pi\alpha}2-a^2-ab-b^2
\ge0.
\end{align}
\end{thm}

Let $p=a+b$ and $q=ab.$
Then, under the assumption of Theorem \ref{thm:ss},
$a$ and $b$ are real numbers precisely if $(a-b)^2=p^2-4q\ge0.$
Otherwise, $a=\bar b=s+it$ for some $s,t\in\R$
and the following result follows from the last theorem.

\begin{cor}\label{cor:ss}
Let $1/3<\alpha<1$ and $s,t\in\R.$
Then the function
$$
f(z)=z\Gauss(s+it,s-it; 2s+1;z),\quad z\in\D,
$$
is strongly starlike of order $\alpha$ if $(s,t)$ is contained
in the closed ellipse given by
$$
\left(s-2\sin^2\frac{\pi\alpha}2\right)^2
+\frac13\left(4\sin^2\frac{\pi\alpha}2-1\right)t^2
\le\sin^2\frac{\pi\alpha}2\left(4\sin^2\frac{\pi\alpha}2-1\right).
$$
\end{cor}

In the next section, we prove Theorem \ref{thm:nec}.
Section 3 will be devoted to proofs of the other results in this section.
We will give some more corollaries and examples in the final section.

\section{Proof of Theorem \ref{thm:nec}}

For the proof, we recall a couple of important formulae of hypergeometric
functions.
For details, the reader can consult monographs \cite{Temme:sp} by Temme
and \cite{WW:anal} by Whittaker and Watson.
As is well known, the hypergeometric function $F(z)=\Gauss(a,b;c;z)$
is characterized as the solution to the hypergeometric differential equation
\begin{equation}\label{eq:Gauss}
(1-z)zF''(z)+[c-(a+b+1)z]F'(z)-abF(z)=0
\end{equation}
with the initial condition $F(0)=1.$
We also note the following relation which readily follows from
the form of the hypergeometric series:
\begin{equation}\label{eq:der}
\frac{d}{dz}\Gauss(a,b;c;z)=\frac{ab}{c}\Gauss(a+1,b+1;c+1;z).
\end{equation}
The following formula for $a,b,c\in\C$ with $a+b\ne c$ and
$c\ne0,-1,-2,\dots$
is useful in what follows:
\begin{align}
&\Gauss(a,b;c;z)
\label{eq:G2}
=\frac{\Gamma(c)\Gamma(c-a-b)}{\Gamma(c-a)\Gamma(c-b)}
\Gauss(a,b;a+b-c+1;1-z) \\
&\qquad +(1-z)^{c-a-b}
\frac{\Gamma(c)\Gamma(a+b-c)}{\Gamma(a)\Gamma(b)}
\Gauss(c-a,c-b;c-a-b+1;1-z).
\notag
\end{align}
When $\Re(c-a-b)>0,$ by \eqref{eq:G2}, we see that
the limit of $\Gauss(a,b;c;z)$ exists as $z\to1$
in $\D$ and evaluated as
\begin{equation}\label{eq:Euler}
\Gauss(a,b;c;1)=\frac{\Gamma(c)\Gamma(c-a-b)}{\Gamma(c-a)\Gamma(c-b)}.
\end{equation}

When $\Re(c-a-b)<0,$ the asymptotic behaviour of $\Gauss(a,b;c;z)$
can be understood via the expression \eqref{eq:G2};
namely, if $\Re(c-a-b)<0,$
\begin{equation}\label{eq:as}
\Gauss(a,b;c;z)=
\frac{\Gamma(c)\Gamma(a+b-c)}{\Gamma(a)\Gamma(b)}(1-z)^{c-a-b}
+O(|1-z|^{\varepsilon})
\end{equation}
as $z\to1$ in $\D,$ where $\varepsilon=\min\{\Re(c-a-b)+1,0\}.$
In the zero-balanced case when $a+b=c,$ we have the following asymptotic
formula due to Ramanujan:
\begin{equation}\label{eq:Ram}
\Gauss(a,b;a+b;z)=\frac{\Gamma(c)}{\Gamma(a)\Gamma(b)}\big(R(a,b)-\log(1-z)\big)
+O\left(|1-z|\log\frac1{|1-z|}\right)
\end{equation}
as $z\to1$ in $\D,$ where
$$
R(a,b)=2\psi(1)-\psi(a)-\psi(b)
$$
and $\psi(x)=\Gamma'(x)/\Gamma(x)$ denotes the digamma function.

We denote by $\D(a,r)$ the open disk $|z-a|<r.$
The next result describes the cluster set
$$
C_1(F)=\bigcap_{0<\delta<1}\overline{F(\D\cap\D(1,\delta))}
$$
of $F(z)=\Gauss(a,b;c;z)$ in the
point $z=1$ in the remaining case when $\Re(c-a-b)=0$ and $c-a-b\ne0.$
We note the simple fact that $C_1(\varphi f)=C_1(f)$ for an analytic
function $f$ on $\D$ whenever $\varphi$ is analytic on $\D$ and
has (unrestricted) limit $1$ as $z\to1$ in $\D.$

\begin{lem}\label{lem:cluster}
Assume that $c=a+b+is$ for an $s\in\R\setminus\{0\}.$
For the function $F(z)=\Gauss(a,b;c;z),$
the cluster set in the point $z=1$ is given by
$$
C_1(F)=\{w\in\C: Re^{-\pi |s|/2}\le |w-w_0|\le Re^{\pi|s|/2}\},
$$
where
$$
w_0
=\frac{\Gamma(c)\Gamma(c-a-b)}{\Gamma(c-a)\Gamma(c-b)}
=\frac{\Gamma(a+b+is)\Gamma(is)}{\Gamma(a+is)\Gamma(b+is)}
$$
and
$$
R=\left|\frac{\Gamma(c)\Gamma(a+b-c)}{\Gamma(a)\Gamma(b)}\right|
=\left|\frac{\Gamma(a+b+is)\Gamma(is)}{\Gamma(a)\Gamma(b)}\right|.
$$
\end{lem}

\begin{pf}
In view of the formula \eqref{eq:G2}, it is enough to look at the function
$$
g(z)=(1-z)^{c-a-b}=(1-z)^{is}=\exp\big[-s\arg(1-z)+i s\log|1-z|\big].
$$
Observe that $\arg g(z)=s\log|1-z|$ is unbounded when $z\to1$
whereas $-\pi/2<\arg(1-z)<\pi/2.$
It is thus easy to deduce the relation
$C_1(g)=\{w: e^{-\pi|s|/2}\le |w|\le e^{\pi|s|/2}\}.$
Since $\overline{\Gamma(is)}=\Gamma(-is),$
the required assertion now follows.
\end{pf}

In particular, $\Gauss(a,b;c;z)$ is bounded on $\D$ in the last case.
As for boundedness of $\Gauss(a,b;c;z),$ we can summarize the above observations.

\begin{lem}\label{lem:bdd}
The hypergeometric function $\Gauss(a,b;c;z)$ is bounded on the unit disk
$\D$ precisely when $\Re(c-a-b)\ge0$ and $c-a-b\ne0.$
\end{lem}

\begin{pf}[Proof of Theorem \ref{thm:nec}]
We put $F(z)=\Gauss(a,b;c;z)$ and $f(z)=zF(z).$
Then
$$
h(z):=\frac{zf'(z)}{f(z)}=1+\frac{zF'(z)}{F(z)}
=1+\frac{abz}c\cdot\frac{\Gauss(a+1,b+1;c+1;z)}{\Gauss(a,b;c;z)}.
$$

Case (i):
When $\Re(c-a-b)>1,$ by \eqref{eq:Euler} for $F$ and $F',$
we compute
$$
h(1)=1+\frac{ab}{c}\cdot
\frac{\Gamma(c+1)\Gamma(c-a-b-1)}{\Gamma(c-a)\Gamma(c-b)}\cdot
\frac{\Gamma(c-a)\Gamma(c-b)}{\Gamma(c)\Gamma(c-a-b)}
=1+\frac{ab}{c-a-b-1},
$$
where we have used the fundamental relation $\Gamma(x+1)=x\Gamma(x)$
for the Euler gamma function.
Since $\Re[e^{-i\lambda}h(1)]\ge0$ by assumption, we have the inequality in (i).

Case (ii):
We next assume that $c-a-b=1+is$ for $s\in\R\setminus\{0\}.$
Note that $F(z)$ has a finite limit as $z\to1$ in $\D.$
Applying Lemma \ref{lem:cluster} to $\Gauss(a+1,b+1;c+1;z),$ we see that
the cluster set $C_1(h)$ is the closed annulus
$Re^{-\pi|s|/2}\le |w-w_1|\le Re^{\pi|s|/2},$ where
$$
w_1=1+\frac{ab}{c F(1)}\cdot
\frac{\Gamma(c+1)\Gamma(c-a-b-1)}{\Gamma(c-a)\Gamma(c-b)}
=1+\frac{ab}{is}
$$
and
$$
R=\left|\frac{ab}{c F(1)}\cdot
\frac{\Gamma(c+1)\Gamma(a+b-c+1)}{\Gamma(a+1)\Gamma(b+1)}\right|
=\left|\frac{\Gamma(-is)\Gamma(c-a)\Gamma(c-b)}{is\Gamma(is)\Gamma(a)\Gamma(b)}
\right|.
$$
Since $|\Gamma(-is)|=|\Gamma(is)|,$ we see that $Re^{\pi|s|/2}$
coincides $R_1$ given in the assertion.
Since the annulus $Re^{-\pi|s|/2}\le |w-w_1|\le Re^{\pi|s|/2}=R_1$
is contained in the closed half-plane
$\Re(e^{-i\lambda}w)\ge0,$ we get the inequalities in the assertion.

Case (iii):
Assume that $c-a-b=1.$
By \eqref{eq:Ram}, we have the asymptotic formula
$$
h(z)=1+\frac{ab}{cF(1)}\cdot\frac{\Gamma(c+1)}{\Gamma(a+1)\Gamma(b+1)}(-\log(1-z)+O(1))
=ab\log\frac1{1-z}+O(1)
$$
as $z\to1$ in $\D.$
Thus we see that the condition in the assertion is necessary.

Case (iv):
Put $c-a-b=\alpha+i\beta$ with $0\le\alpha<1$ and $\beta\in\R.$
We first assume that $0<\alpha<1.$
By \eqref{eq:Euler} and \eqref{eq:as}, we obtain
$$
h(z)=1+A(1-z)^{c-a-b-1}(1+o(1))
=1+A(1-z)^{\alpha-1+i\beta}(1+o(1))
$$
as $z\to1$ in $\D,$ where
$$
A=\frac{ab}{cF(1)}\cdot
\frac{\Gamma(c+1)\Gamma(a+b-c+1)}{\Gamma(a+1)\Gamma(b+1)}
=\frac{\Gamma(c-a)\Gamma(c-b)\Gamma(a+b+1-c)}{\Gamma(a)\Gamma(b)\Gamma(c-a-b)}.
$$
Since $|\arg(1-z)|<\pi/2$ in $|z|<1,$
we see that
$$
|h(z)|=|A|\exp\left\{-(1-\alpha)\log|1-z|-\beta\arg(1-z)\right\}(1+o(1))
\to+\infty
$$
as $z\to1$ and that
$$
\arg h(z)=\beta\log|1-z|+(\alpha-1)\arg(1-z)+\arg A
\to-\sgn(\beta)\infty
$$
as $z\to1$ if $\beta\ne0.$
Hence, the image $h(\D)$ cannot be contained in the half-plane
$\Re(e^{-i\lambda}w)>0$ if $\beta\ne0.$
Therefore, the assumption that $f$ is $\lambda$-spirallike
implies that $\beta=0.$
Moreover, we should have $\lambda-\pi/2\le (\alpha-1)\pi/2+\arg A
\le (1-\alpha)\pi/2+\arg A\le\lambda+\pi/2,$
which implies the inequality in the assertion,
where one should use the fact that $\arg\Gamma(a+b+1-c)=\arg\Gamma(c-a-b)=0.$

Secondly, we assume $\alpha=0$ and $\beta\ne0.$
Then, by \eqref{eq:G2},
$$
h(z)
=\frac{A+o(1)}{1-z}\cdot\frac{(1-z)^{i\beta}}{B+(1-z)^{i\beta}}
=\frac{A+o(1)}{1-z}\cdot\frac{1}{B+(1-z)^{-i\beta}}
$$
as $z\to1,$ where $A$ and $B$ are nonzero complex numbers.
As we saw above, the cluster set of $(1-z)^{-i\beta}$ in the point $z=1$
is the annulus $e^{-\pi|\beta|/2}\le |w|\le e^{\pi|\beta|/2}.$
Therefore, the inequality $e^{\pi|\beta|/2}\le|B|$ should hold.
For a fixed $\theta\in(-\pi/2,\pi/2),$ consider the curve
$z_\theta(t)=1-te^{i\theta}.$
Note that $z_\theta(t)\in\D$ for $t\in(0,t_\theta)$ for a positive number
$t_\theta.$
We now have
$$
h(z_\theta(t))
=\frac{A+o(1)}{te^{i\theta}}\cdot
\frac{1}{B+e^{\beta\theta}\exp[i\beta\log(1/t)]}\quad(t\to0^+).
$$
Therefore, if we denote by $\Phi_+(\theta)$ and $\Phi_-(\theta)$
the upper and lower limits of $\arg h(z_\theta(t))$ as $t\to0^+$ respectively,
we obtain
$$
\Phi_\pm(\theta)=\arg A-\theta-\arg B\pm\arcsin\frac{e^{\beta\theta}}{|B|}.
$$
Therefore,
$$
\lim_{\theta\to\pi/2^-}\big[\Phi_+(-\theta)-\Phi_-(\theta)\big]
=\pi+\arcsin\frac{e^{-\pi\beta/2}}{|B|}+\arcsin\frac{e^{\pi\beta/2}}{|B|}>\pi,
$$
which implies that the image $h(\D)$ cannot be contained in the half-plane
$\Re(e^{-i\lambda}w)>0.$

Finally, we assume that $\alpha=\beta=0.$
By \eqref{eq:as} and \eqref{eq:Ram}, we obtain
$$
h(z)=\frac{1+o(1)}{-(1-z)\log(1-z)}
$$
as $z\to1$ in $\D.$
By using the above curve $z_\theta(t),$ we compute
$h(z_\theta(t))=(1+o(1))e^{-i\theta}t\inv/$
$[-\log(1/t)-i\theta].$
In particular, $\arg h(z_\theta(t))\to -\theta$ as $t\to0^+$
for $-\pi/2<\theta<\pi/2.$
Therefore, the image $h(\D)$ can be contained only in the half-plane
$\Re w>0;$ in other words, $\lambda=0.$
Thus the assertion is deduced in this case, too.

Case (v):
Assume that $\Re(c-a-b)<0.$
By \eqref{eq:as}, we have
$$
h(z)=\frac{a+b-c+o(1)}{1-z}
$$
as $z\to1$ in $\D.$
Thus $\lambda=\arg(a+b-c)$ as required in order that the image $h(\D)$ is
contained in the half-plane $\Re(e^{-i\lambda}w)>0.$
\end{pf}

\section{Proofs of the other results}

As in a paper \cite{Lewis79} of Lewis (see also the proof of
\cite[Theorem 2.12]{Rus:conv}), our proof will be based on
the following lemma due to Jack \cite{Jack71} and the hypergeometric
differential equation \eqref{eq:Gauss}.
The strategy and computations are largely overlapped with those
in K\"ustner \cite[\S 3]{Kus07} at least when $\lambda=0$
(and his paper, more generally, deals with the estimate of $\sigma(f)$).
We note, however, that final conclusions are given in \cite{Kus07}
only for real parameters $a,b,c.$

\begin{lem}[Jack's lemma]\label{lem:Jack}
Let $\omega$ be a non-constant analytic function on $\D$ with $\omega(0)=0$.
If the maximum value of $|\omega(z)|$ on the circle $|z|=r$ with $0<r<1$ is
attained at a point $z_0$ on the circle,
then $z_0\omega'(z_0)=k\omega(z_0)$ for some $k\geq 1$.
\end{lem}

Let $F(z)=\Gauss(a,b;c;z)$ and $f(z)=zF(z).$
For a while, without any assumptions on the parameters $a,b,c,$
we try to show that $f$ is $\lambda$-spirallike for a fixed
$\lambda\in(-\pi/2,\pi/2).$
We define a meromorphic function $p$ on $\D$ with $p(0)=1$
by the relation
\begin{equation*}
\frac{zf'(z)}{f(z)}
=e^{i\lambda}\big(p(z)\cos\lambda -i\sin\lambda\big).
\end{equation*}
In view of the formula $zf'(z)/f(z)=1+zF'(z)/F(z),$ we obtain
$$
\frac{zF'(z)}{F(z)}=e^{i\lambda}\big(p(z)\cos\lambda -i\sin\lambda\big)-1
=\mu(p(z)-1),
$$
where
$$
\mu=e^{i\lambda}\cos\lambda.
$$
Differentiating the derived formula $zF'(z)=\mu(p(z)-1)F(z),$ we get
$$
zF''(z)+F'(z)=\mu p'(z)F(z)+\mu(p(z)-1)F'(z).
$$
In conjunction with \eqref{eq:Gauss}, we have
$$
[1-c+(a+b)z-\mu(1-z)(p(z)-1)]F'(z)=[\mu(1-z)p'(z)-ab]F(z),
$$
which further leads to
\begin{equation}\label{eq:eq}
[1-c+(a+b)z-\mu(1-z)(p(z)-1)]\mu(p(z)-1)=\mu(1-z)zp'(z)-abz.
\end{equation}

To verify $\lambda$-spirallikeness of $f(z)$,
we need to show that $\Re p(z)>0$ for $z\in\D.$
It is equivalent to the condition that the meromorphic function
\begin{equation}\label{eq:omega}
\omega(z)=\frac{p(z)-1}{p(z)+1}
\end{equation}
satisfies $|\omega|<1$ on $\D.$
Suppose, to the contrary, that there exists
a $z_0\in\D$ such that $|\omega(z_0)|=1$ and that $|\omega(z)|<1$ for
$|z|<r_0:=|z_0|.$
Then Lemma \ref{lem:Jack} implies that
$z_0\omega'(z_0)/\omega(z_0)=k\geq 1.$

If $\omega(z_0)=1,$ then the meromorphic function 
$q=1/p=(1-\omega)/(1+\omega)$ satisfies
$q(z_0)=0$ and $z_0q'(z_0)=-z_0\omega'(z_0)/2=-k/2.$
We substitute $p=1/q$ into \eqref{eq:eq} to obtain the relation
$$
\big[q(z)\{1-c+(a+b)z\}-\mu(1-z)(1-q(z))\big]\mu(1-q(z))
=-\mu(1-z)zq'(z)-abzq(z)^2.
$$
Letting $z=z_0,$ we obtain $k=-2\mu,$ which is impossible.
Thus we have $\omega(z_0)\ne1.$
Hence, the function
\begin{equation*}
p(z)=\frac{1+\omega(z)}{1-\omega(z)}
\end{equation*}
is analytic at $z=z_0$ and satisfies $\Re p(z_0)=0.$
Taking the logarithmic derivative of the both sides of \eqref{eq:omega},
we have
$$
\frac{2z_0p'(z_0)}{p(z_0)^2-1}=\frac{z_0\omega'(z_0)}{\omega(z_0)}=k.
$$
Thus, we can write
\begin{equation}\label{eq:s}
p(z_0)=i s
\aand
z_0p'(z_0)=-k(s^2+1)/2\le -(s^2+1)/2
\end{equation}
for some $s\in\R.$
We put
\begin{equation}\label{eq:st}
\sigma
=\mu(is-1)
\aand
\tau
=\mu k(s^2+1)/2.
\end{equation}
Letting $z=z_0$ in \eqref{eq:eq} and recalling \eqref{eq:s},
we obtain the relation
\begin{equation*}
\tau-\sigma(\sigma+c-1)=[\tau-(\sigma+a)(\sigma+b)]z_0.
\end{equation*}
Therefore, we will get a contradiction if the inequality
\begin{equation}\label{eq:ineq}
\big|\tau-\sigma(\sigma+c-1)\big|\ge \big|\tau-(\sigma+a)(\sigma+b)\big|
\end{equation}
holds and if equalities
\begin{equation*}
\tau-\sigma(\sigma+c-1)=\tau-(\sigma+a)(\sigma+b)=0
\end{equation*}
never hold simultaneously for any $s\in\R$ and $k\ge1,$
where $\sigma, \tau$ are given in \eqref{eq:st}.

Hence, it is enough to show \eqref{eq:ineq} and 
$\sigma(\sigma+c-1)\ne (\sigma+a)(\sigma+b)$
to prove $\lambda$-spirallikeness of the function $f(z)=z\Gauss(a,b;c;z).$

Fix $s$ (and thus $\sigma$) for a while.
The inequality \eqref{eq:ineq} means exactly that the point $\tau$
is contained in the half-plane $H$ bounded by the perpendicular bisector
of the two points $A=\sigma(\sigma+c-1)$ and $B=(\sigma+a)(\sigma+b),$
which contains the point $B,$ provided that $A\ne B.$
Note that the point $\tau=\mu k(s^2+1)/2$ with
$k\ge1$ may vary on the ray emanating from the point $\tau_1=\mu(s^2+1)/2$
with the direction $\mu.$
Hence, the ray is contained in $H$ precisely when $\tau_1\in H$ and
$|\arg(B-A)-\arg\mu|\le\pi/2.$
The second condition and the condition $A\ne B$ follow from the inequality
\begin{equation}\label{eq:direction}
\Re[(B-A)\bar\mu]=\Re\left[|\mu|^2(is-1)(a+b+1-c)+ab\bar\mu\right]>0.
\end{equation}
(We remark that the second condition follows from the weaker
inequality $\Re[(B-A)\bar\mu]\ge0$ when $A\ne B$ is already established
by another way.)
Since the inequality \eqref{eq:direction} should hold for any $s\in\R,$ 
the condition $\Im(a+b+1-c)=0$ is required.
Thus $p:=a+b+1-c$ is a real number.
We should also have the inequality $\Re[ab\bar\mu]> p|\mu|^2.$
In other words, $\Re\big[e^{-i\lambda}ab\big]>p\cos\lambda;$
equivalently,
$$
\Re\big[e^{-i\lambda}(ab-p)\big]>0.
$$

We next consider the first condition $\tau_1\in H;$ namely,
$|\tau_1-A|\ge|\tau_1-B|.$
By squaring, we see that it is equivalent to validity of the inequality
\begin{equation}\label{eq:t0}
|A|^2-|B|^2-2\Re\left[(A-B)\overline{\tau_1}\right]\ge0.
\end{equation}
A substitution of the concrete forms of $A, B, \sigma$ and $\tau_1$ gives us
$$
2\Re\left[(A-B)\overline{\tau_1}\right]
=(s^2+1)(p|\mu|^2-\Re[ab\bar\mu]),
$$
which is a quadratic polynomial in $s.$
On the other hand,
\begin{align*}
|A|^2-|B|^2
&= |\sigma|^2(|\sigma|^2+2\Re(c-1)\bar\sigma +|c-1|^2) \\
&\quad -(|\sigma|^2+2\Re[a\bar\sigma] +|a|^2)(|\sigma|^2+2\Re[b\bar\sigma]+|b|^2)
 \\
&=-2\Re[(a+b+1-c)\bar\sigma]|\sigma|^2+(|c-1|^2-|a|^2-|b|^2)|\sigma|^2 \\
&\quad
-(2\Re[a\bar\sigma]+|a|^2)(2\Re[b\bar\sigma]+|b|^2).
\end{align*}
Since the first term in the last expression is equal to
$$
-2p\Re\bar\sigma|\sigma|^2=2p(s\Im\mu+\Re\mu)|\mu|^2(s^2+1)
=2p(\Im\mu)|\mu|^2s^3+O(s^2)
$$
and the other terms are polynomials in $s$ of degree at most $2,$
we need the condition $p\Im\mu=0$ for the inequality \eqref{eq:t0} to hold
for all $s\in\R.$
Hence, the present approach works only when $\lambda=0$ or $p=0,$
which correspond to Theorems \ref{thm:st} and \ref{thm:spl}, respectively.
We are now ready to prove these theorems.

\begin{pf}[Proof of Theorem \ref{thm:st}]
Here, we assume that $\lambda=0.$
Therefore, we now have $\mu=1$ and $\sigma=-1+is.$
For convenience, we write $a=a_1+ia_2$ and $b=b_1+ib_2.$
Substituting these, the left-hand side of \eqref{eq:t0} can be computed as
\begin{align*}
&\qquad
2p(s^2+1)+(|c-1|^2-|a|^2-|b|^2)(s^2+1) \\
&\qquad -(-2a_1-2sa_2+|a|^2)(-2b_1-2sb_2+|b|^2)+(s^2+1)(\Re[ab]-p) \\
&=(|c-1|^2-|a|^2-|b|^2+p+\Re[ab]-4a_2b_2)s^2
-2(2a_1b_2+2a_2b_1-a_2|b|^2-b_2|a|^2)s \\
&\qquad+(|c-1|^2-|a|^2-|b|^2+p+\Re[ab]-4a_1b_1+2a_1|b|^2+2b_1|a|^2-|a|^2|b|^2)\\
&=(|c-1|^2-|a+b|^2+p+3\Re[ab])s^2-2\Im[\bar a\bar b(a+b-2)]s \\
&\qquad+|c-2|^2-|a-1|^2|b-1|^2-p+\Re[ab].
\end{align*}
Since the above quadratic polynomial in $s$ is non-negative,
the assertion follows.
\end{pf}

\begin{rem}\label{rem:eq}
In Theorem \ref{thm:st}, we assumed the strict inequality $\Re[ab]>p.$
We can, however, weaken the assumption to $\Re[ab]\ge p$
by a limiting argument in some cases.
For instance, we assume that $\Re[ab]=p$ and that
$L, N $ and $LN-M^2$ are all positive.
Then, for $\varepsilon>0,$ we consider the function
$f_\varepsilon(z)=z\Gauss(a,b;c+\varepsilon;z).$
Note that $f_\varepsilon$ converges to the original function $f=f_0$
locally uniformly on $\D$ as $\varepsilon\to0.$
We now observe that $p_\varepsilon=a+b+1-c-\varepsilon$ is real and
$\Re[ab]-p_\varepsilon=\varepsilon>0.$
Moreover, for a sufficiently small $\varepsilon>0,$
the quantities $L_\varepsilon, N_\varepsilon, 
L_\varepsilon N_\varepsilon-M_\varepsilon^2$ corresponding to $f_\varepsilon$
are all still positive.
Therefore, by the theorem, we conclude that $f_\varepsilon$ is starlike.
Since starlikeness is preserved by locally uniform convergence,
we see that $f$ is starlike.
On the other hand, this procedure does not necessary work
when the quadratic form $Ls^2+2Mst+Nt^2$ is degenerate.
See Remark \ref{rem:ex} below.
\end{rem}

\begin{pf}[Proof of Theorem \ref{thm:spl}]
We next complete the proof of Theorem \ref{thm:spl}.
We can harmlessly assume that $ab\ne0.$
By the form of the function, $p=a+b+1-(a+b+1)=0.$
Then we have
$B-A=p\sigma+ab=ab\ne0$ and
\begin{align*}
|A|^2-|B|^2
&=(|a+b|^2-|a|^2-|b|^2)|\sigma|^2
-(2\Re[a\bar\sigma]+|a|^2)(2\Re[b\bar\sigma]+|b|^2) \\
&=2|\sigma|^2\Re[a\bar b]
-(2\Re[a\bar\sigma]+|a|^2)(2\Re[b\bar\sigma]+|b|^2).
\end{align*}
Thus, by recalling $\sigma=(is-1)\mu,$ we compute
the left-hand side of \eqref{eq:t0} as
\begin{align*}
& 2(s^2+1)|\mu|^2\Re[a\bar b]+(s^2+1)\Re[ab\bar\mu] \\
& -(2s\Im[a\bar\mu]-2\Re[a\bar\mu]+|a|^2)
(2s\Im[b\bar\mu]-2\Re[b\bar\mu]+|b|^2) \\
= \quad
& ~(\Re[ab\bar\mu]+2\Re[a\bar\mu\overline{b\bar\mu}]
-4\Im[a\bar\mu]\Im[b\bar\mu])s^2 \\
& -2\big[\Im[a\bar\mu](-2\Re[b\bar\mu]+|b|^2)
+\Im[b\bar\mu](-2\Re[a\bar\mu]+|a|^2)\big]s \\
& +\Re[ab\bar\mu]+2\Re[a\bar\mu\overline{b\bar\mu}]
-(-2\Re[a\bar\mu]+|a|^2)(-2\Re[b\bar\mu]+|b|^2) \\
= \quad
& ~(\Re[ab\bar\mu]+2\Re[ab\bar\mu^2])s^2
-2\Im\big[-2ab\bar\mu^2+ab\bar b\bar\mu+ba\bar a\bar\mu\big]s \\
& +\Re[ab\bar\mu]-2\Re[ab\bar\mu^2]+2\Re[ab(\bar a+\bar b)\bar\mu]-|a|^2|b|^2 \\
= \quad
&\cos\lambda(Ls^2-2Ms+N),
\end{align*}
where $L, M, N$ are given in the assertion of Theorem \ref{thm:spl}.
In the above, we used the relation $2\mu=1+e^{2i\lambda}.$

Finally, we observe that the assumptions in the theorem imply that
$c=a+b+1\ne0,-1,-2.\dots.$
Suppose, to the contrary, that $a+b=-k$ for some $k\in\{1,2,3,\dots\}.$
Then
$$
N=-2k\Re[e^{-i\lambda}ab]-\Re[e^{-3i\lambda}ab]-|ab|/\cos\lambda
\le -2k\Re[e^{-i\lambda}ab].
$$
Since the last term is negative by condition (i), we have a contradiction.
Hence, we conclude that $a+b\ne-1,-2,\dots.$
The proof is now complete.
\end{pf}

In order to prove Theorem \ref{thm:convex}, it is enough to note
the following fact which follows from Alexander's theorem
(see \cite[Theorem 2.12]{Duren:univ}).

\begin{lem}\label{lem:Alex}
Let $a,b,c$ be complex numbers with $ab\ne0,~c\ne0,-1,-2,\dots.$
The Gaussian hypergeometric function $g(z)=\Gauss(a,b;c;z)$ is convex
if and only if $f(z)=z\Gauss(a+1,b+1;c+1;z)$ is starlike.
\end{lem}

\begin{pf}
Note the relation  $f(z)=(ab/c)zg'(z)$ by the formula \eqref{eq:der}.
By taking the logarithmic derivatives of the both sides,
we obtain the relation
$$
\frac{zf'(z)}{f(z)}=1+\frac{zg''(z)}{g'(z)},
$$
from which the assertion follows.
\end{pf}

\begin{pf}[Proof of Corollary \ref{cor:1}]
In order to apply Theorem \ref{thm:spl},
under the assumptions of the corollary, we compute
\begin{align*}
L&=q\Re\big[(2+e^{-2i\lambda})\big]=q(2+\cos2\lambda), \\
M&=q\Im\big[\bar a+\bar b-2e^{-i\lambda}\cos\lambda\big]
=-q(\Im[a+b]-2\sin\lambda\cos\lambda),
\aand \\
N&=q\Re\big[2\bar a+2\bar b-e^{-2i\lambda}
-e^{i\lambda}\bar a\bar b/\cos\lambda\big]
=q(2\Re[a+b]-\cos2\lambda-q/\cos\lambda).
\end{align*}
Thus the assertion follows.
\end{pf}

\begin{pf}[Proof of Corollary \ref{cor:2}]
Similarly, assuming that $q=ab$ is positive, we compute
\begin{align*}
L&=q\Re\big[e^{-i\lambda}(2+e^{-2i\lambda})\big]
=q\cos\lambda(4\cos^2\lambda-1), \\
M&=q\Im\big[e^{-i\lambda}(\bar a+\bar b-2e^{-i\lambda}\cos\lambda)\big]
=-q(\Im[e^{i\lambda}(a+b)]-2\sin2\lambda\cos\lambda),
\aand \\
N&=q\Re\big[e^{-i\lambda}(2\bar a+2\bar b-e^{-2i\lambda}
-e^{i\lambda}\bar a\bar b/\cos\lambda)\big]
=q(2\Re[e^{i\lambda}(a+b)]\cos\lambda-\cos3\lambda-q).
\end{align*}
The assertion now follows from Theorem \ref{thm:spl}.
\end{pf}

By using the last corollary, we are able to show Theorem \ref{thm:ss}.

\begin{pf}[Proof of Theorem \ref{thm:ss}]
Let $\lambda=(1-\alpha)\pi/2~(<\pi/3).$
When $p=a+b\in\R$ and $q=ab>0,$ the condition \eqref{eq:spl} reads
\begin{align*}
&~(4\cos^2\lambda-1)(2p\cos^2\lambda-4\cos^4\lambda
+3\cos^2\lambda-q)-(p\sin\lambda-2\sin2\lambda\cos\lambda)^2 \\
=&~
-p^2\sin^2\lambda+6p\cos^2\lambda+q(1-4\cos^2\lambda)-3\cos^2\lambda
\ge0,
\end{align*}
which is nothing but \eqref{eq:ss}.
Noting the relation $\cos\lambda=\sin[\pi\alpha/2],$ we see that \eqref{eq:spl}
is equivalent to \eqref{eq:ss}.
Since the condition \eqref{eq:ss} is unchanged if we replace $\lambda$
by $-\lambda,$ we conclude that the function $z\Gauss(a,b;a+b+1;z)$
is contained in the class $\spl(\lambda)\cap\spl(-\lambda)=\st(\alpha).$
\end{pf}

\section{Some examples}
This section is devoted to giving examples of spirallike shifted
hypergeometric functions.

Let $a=2$ in Theorem \ref{thm:st}.
Then $\Re[ab]-p=2\Re b-(b-c+3)=\bar b+c-3$ and
$L, M, N$ have the simple forms $L=N=\Re[c-b]\cdot\Re[b+c-3],~M=0$
so that $LN-M^2=L^2\ge0.$
Therefore the next result follows from Theorem \ref{thm:st}.

\begin{cor}
Let $b,c,s$ be real numbers with $3\le b+c$ and $b\le c.$
Then the function $f(z)=z\Gauss(2,b+is;c+is;z)$ is starlike and
the function $g(z)=z\Gauss(1,b+is;c+is;z)$ is convex.
\end{cor}

\begin{pf}
As is accounted above, starlikeness of $f$ follows from Theorem \ref{thm:st}
when $3<b+c.$
When $b+c=3,$ we first apply the theorem to the function
$z\Gauss(2,b+is;c+\varepsilon+is;z)$ for $\varepsilon>0$ and let $\varepsilon
\to0.$
The second assertion follows from the Alexander relation $zg'(z)=f(z)$
(see Lemma \ref{lem:Alex}) which can be checked by comparing the
coefficients of the power series expansions.
\end{pf}

\begin{rem}\label{rem:ex}
It is noteworthy here that the triple $(a,b,c)=(2,b,3-\bar b)$
satisfies $\Re[ab]=p, L=M=N=0.$
Therefore, as far as $b\ne 3,4,5,\dots,$ the function
$f(z)=z\Gauss(2,b;3-\bar b;z)$ satisfies all the assumptions in
Theorem \ref{thm:st} with (ii) $\Re[ab]>p$ being replaced by $\Re[ab]=p.$
On the other hand, it is necessary for $f$ to be starlike
that $|f''(0)/2|=|2b/(3-\bar b)|\le 2,$ which is equivalent to
$\Re b\le 3/2.$
Therefore, the above function $f$ with $\Re b>3/2$ tells us that
we cannot replace the condition (ii) by $\Re[ab]\ge p$ in general.
\end{rem}

Put $\gamma=b-1+is.$
When $c=b+1,$ the function $g(z)$ in the last corollary takes the form
$$
z\Gauss(1,\gamma+1;\gamma+2;z)
=\sum_{j=1}^{\infty}\frac{\gamma+1}{\gamma+j}z^j.
$$
The corollary implies that it is convex
if $\Re\gamma\geq 0.$
Note that this is contained in Theorem 5 with $n=1$ in Ruscheweyh \cite{Rus75}.

\medskip

Similarly, let $a=2$ in Theorem \ref{thm:spl}
and $b=Re^{i\lambda}$ with $R>0.$
Then $L, M, N$ in Theorem \ref{thm:spl} are expressed by
\begin{align*}
L&=2R(2+\cos2\lambda)=2R(1+2\cos^2\lambda), \\
M&=-2R(R-2\cos\lambda)\sin\lambda,\\
N&=2R(2R\cos\lambda+5-2\cos^2\lambda-2R/\cos\lambda).
\end{align*}
Therefore,
\begin{align*}
&\quad~ LN-M^2~
=4R^2\left(
-R^2\sin^2\lambda-2R\sin\lambda\tan\lambda+5-4\cos^2\lambda
\right) \\
&=-\frac{4R^2(R\sin\lambda\cos\lambda+\sin\lambda-3+2\sin^2\lambda)%
(R\sin\lambda\cos\lambda+\sin\lambda+3-2\sin^2\lambda)}{\cos^2\lambda}.
\end{align*}
For simplicity, we may assume that $\lambda>0.$
(Note that $f(z)$ is $\lambda$-spirallike if and only if
$\overline{f(\bar z)}$ is $-\lambda$-spirallike.)
Then $LN-M^2\ge0$ if and only if
$R\sin\lambda\cos\lambda\le -\sin\lambda+3-2\sin^2\lambda
=(1-\sin\lambda)(3+2\sin\lambda).$
Simpler conditions for spirallikeness can now be obtained by
Theorem \ref{thm:spl} as follows.

\begin{cor}\label{cor:spl}
Let $0<\lambda<\pi/2.$
Then the function $z\Gauss(2,b;3+b;z)$ with $b=Re^{i\lambda}$
is $\lambda$-spirallike if
$$
0<R\le\frac{(1-\sin\lambda)(3+2\sin\lambda)}{\sin\lambda\cos\lambda}.
$$
\end{cor}

We next let $a=2e^{i\lambda}\cos\lambda$ and $b>0$ in Theorem \ref{thm:spl}
to obtain the following.

\begin{cor}
Let $\lambda$ be a real number with $0<|\lambda|<\pi/2.$
Then the shifted hypergeometric function 
$z\Gauss(2e^{i\lambda}\cos\lambda,b; 2e^{i\lambda}\cos\lambda+b+1;z)$
is $\lambda$-spirallike for any constant $b>0.$
\end{cor}

Indeed, in this case, we check condition (iii) in the theorem by
$L=2b(2+\cos2\lambda)\cos\lambda>0, M=0,
N=2b(2\cos^2\lambda+1)\cos\lambda>0.$

Next, let $a+b=s\in\R$ and $ab=qe^{i\lambda}$ with $q>0.$
Then we see the following by applying Corollary \ref{cor:1}.

\begin{cor}\label{cor:s}
Let $\lambda\in(-\pi/2,\pi/2)$ with $\lambda\ne0.$
Suppose that $a,b\in\C$ satisfy $a+b=s\in\R$ and $ab=qe^{i\lambda}$
for some $q>0.$
Then the function $z\Gauss(a,b;s+1;z)$ is $\lambda$-spirallike if
$$
2s-\frac{q}{\cos\lambda}\ge\frac{1+2\cos2\lambda}{2+\cos2\lambda}
=\frac{4\cos^2\lambda-1}{2\cos^2\lambda+1}.
$$
\end{cor}

\begin{example}
Let $\lambda=\pi/4$ in Corollary \ref{cor:s}.
Then the required inequality takes the form $2s-\sqrt{2}q\ge 1/2.$
For instance, put $a=5(1+2i)/8$ and $b=5(3-i)/4.$
Then $s=a+b=35/8$ and $q=abe^{-\pi i/4}=125(1+i)e^{-\pi i/4}/32=125/16\sqrt2$
satisfy the inequality.
See Figure 1, generated by Mathematica ver.\,10,
for the image of $\D$ under the mapping
$f(z)=z\Gauss(a,b;a+b+1;z)$ in this case.
As the picture suggests us, $f$ is not starlike.
Indeed, the real part of $zf'(z)/f(z)$ assumes the value $-0.0374$
approximately at $z=e^{\pi i/4}.$
\end{example}

\begin{figure}\label{fig:spl}
\begin{center}
\includegraphics[height=.6\textheight]{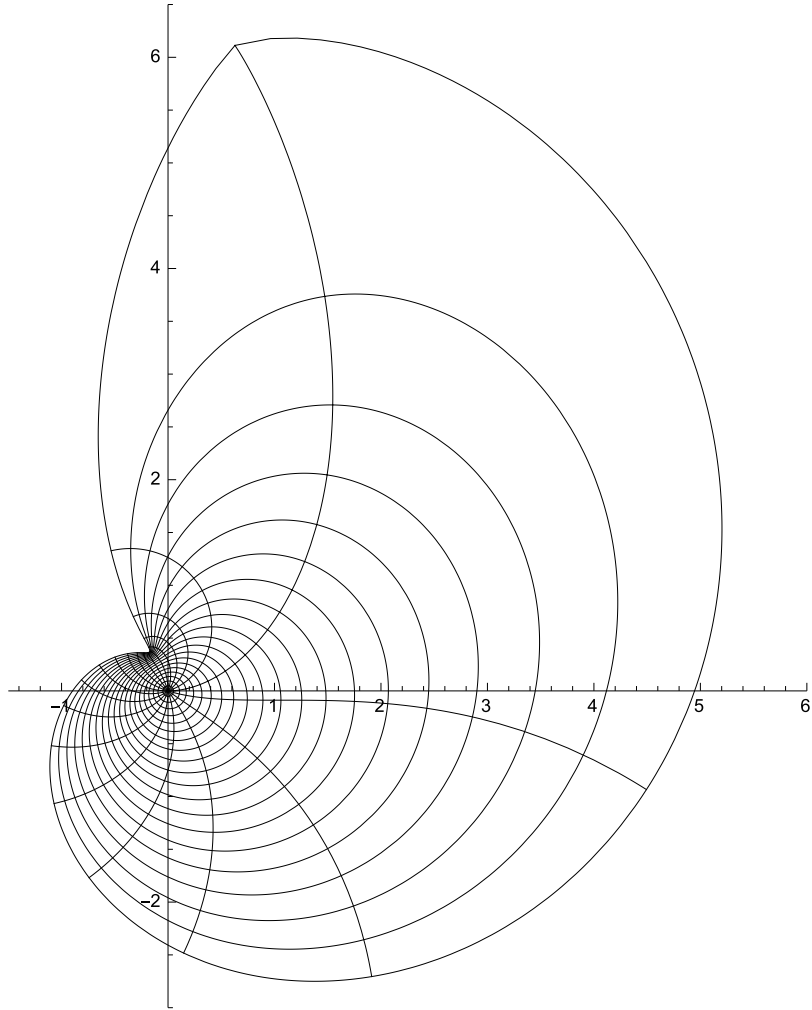}
\caption{Image domain of $z\Gauss(a,b;a+b+1;z)$ with
$a+b=35/8$ and $ab=125e^{\pi i/4}/16\sqrt2$}
\end{center}
\end{figure}

\def\cprime{$'$} \def\cprime{$'$} \def\cprime{$'$}
\providecommand{\bysame}{\leavevmode\hbox to3em{\hrulefill}\thinspace}
\providecommand{\MR}{\relax\ifhmode\unskip\space\fi MR }
\providecommand{\MRhref}[2]{%
  \href{http://www.ams.org/mathscinet-getitem?mr=#1}{#2}
}
\providecommand{\href}[2]{#2}


\begin{thebibliography}{10}

\bibitem{AS:hand}
M.~Abramowitz and I.~A. Stegun, \emph{Handbook of {M}athematical {F}unctions},
  Dover, 1972.

\bibitem{Duren:univ}
P.~L. Duren, \emph{Univalent {F}unctions}, Springer-Verlag, 1983.

\bibitem{HPV10}
P.~H\"asto, S.~Ponnusamy, and M.~Vuorinen, \emph{Starlikeness of the {G}aussian
  hypergeometric functions}, Complex Var. Elliptic Equ. \textbf{55} (2010),
  173--184.

\bibitem{Jack71}
I.~S. Jack, \emph{Functions starlike and convex of order alpha}, J. London
  Math. Soc. \textbf{3} (1971), 469--474.

\bibitem{KS12spl}
Y.~C. Kim and T.~Sugawa, \emph{Correspondence between spirallike functions and
  starlike functions}, Math. Nachr. \textbf{285} (2012), 322--331.

\bibitem{Kus02}
R.~K\"ustner, \emph{Mapping properties of hypergeometric functions and
  convolutions of starlike or convex functions of order $\alpha$}, Comput.
  Methods Funct. Theory \textbf{2} (2002), 597--610.

\bibitem{Kus07}
R.~K{\"u}stner, \emph{On the order of starlikeness of the shifted {G}auss
  hypergeometric function}, J. Math. Anal. Appl. \textbf{334} (2007),
  1363--1385.

\bibitem{Lewis79}
J.~Lewis, \emph{Applications of a convolution theorem to {J}acobi polynomails},
  SIAM J. Math. Anal. \textbf{10} (1979), 1110--1120.

\bibitem{MS61}
E.~P. Merkes and W.~T. Scott, \emph{Starlike hypergeometric functions}, Proc.
  Amer. Math. Soc. \textbf{12} (1961), 885--888.

\bibitem{MM90U}
S.~Miller and P.~Mocanu, \emph{Univalence of {G}aussian and confluent
  hypergeometric functions}, Proc. Amer. Math. Soc. \textbf{119} (1990),
  333--342.

\bibitem{PV01}
S.~Ponnusamy and M.~Vuorinen, \emph{Univalence and convexity properties for
  {G}aussian hypergeometric functions}, Rocky Mountain J. Math. \textbf{31}
  (2001), 327--353.

\bibitem{Rus75}
{St}. Ruscheweyh, \emph{New criteria for univalent functions}, Proc. Amer.
  Math. Soc. \textbf{49} (1975), 109--115.

\bibitem{Rus:conv}
\bysame, \emph{Convolutions in {G}eometric {F}unction {T}heory}, S\'eminaire de
  Math\'ematiques Sup\'erieures, vol.~83, Les Presses de l'Universit\'e de
  Montr\'eal, Montr\'eal, 1982.

\bibitem{RS86}
{St}. Ruscheweyh and V.~Singh, \emph{On the order of starlikeness of
  hypergeometric functions}, J. Math. Anal. Appl. \textbf{113} (1986), 1--11.

\bibitem{Sil93}
H.~Silverman, \emph{Starlike and convexity properties for hypergeometric
  functions}, J. Math. Anal. Appl. \textbf{172} (1993), 574--581.

\bibitem{SugawaDual}
T.~Sugawa, \emph{A self-duality of strong starlikeness}, Kodai Math J.
  \textbf{28} (2005), 382--389.

\bibitem{Temme:sp}
N.~M. Temme, \emph{Special {F}unctions}, John Wiley \& Sons, Inc., New York,
  1996, An introduction to the classical functions of mathematical physics.

\bibitem{WW:anal}
E.~T. Whittaker and G.~N. Watson, \emph{A course of {M}odern {A}nalysis. {A}n
  {I}ntroduction to the {G}eneral {T}heory of {I}nfinite {P}rocesses and of
  {A}nalytic {F}unctions; with an {A}ccount of the {P}rincipal {T}ranscendental
  {F}unctions}, Fourth edition. Reprinted, Cambridge University Press, New
  York, 1962.

\end{thebibliography}

\end{document}